\documentclass[11pt]{amsart}

\usepackage{amsmath,xspace,amssymb,mathrsfs}
\usepackage{color}

\input xy
\xyoption{all}
\xyoption{2cell}
\UseAllTwocells
\CompileMatrices

\newcommand{\Spec}{\operatorname{Spec}}
\renewcommand{\phi}{\varphi}

\newcommand{\Ker}{\operatorname{Ker}}

\newcommand{\Ima}{\operatorname{Im}}

\newcommand{\Max}{\operatorname{Max}}

\newcommand{\Min}{\operatorname{Min}}
\newcommand{\Ann}{\operatorname{Ann}}

\newcommand{\Supp}{\operatorname{Supp}}

\newcommand{\D}{\operatorname{D}}
\newcommand{\Clop}{\operatorname{Clop}}
\newcommand{\Fin}{\operatorname{Fin}}

\newcommand{\Top}{\operatorname{Top}}

\newcommand{\Set}{\operatorname{Set}}
\newcommand{\Ring}{\operatorname{Ring}}
\newcommand{\Mor}{\operatorname{Mor}}
\newcommand{\Su}{\operatorname{S}}

\newtheorem{proposition}{Proposition}[section]
\newtheorem{lemma}[proposition]{Lemma}

\newtheorem{corollary}[proposition]{Corollary}
\newtheorem{theorem}[proposition]{Theorem}

\theoremstyle{definition}
\newtheorem{definition}[proposition]{Definition}
\newtheorem{example}[proposition]{Example}

\newtheorem{remark}[proposition]{Remark}

\usepackage{etoolbox}
\makeatletter
\patchcmd{\@settitle}{\uppercasenonmath\@title}{}{}{}
\patchcmd{\@setauthors}{\MakeUppercase}{}{}{}
\makeatother

\begin{document}

\title{Minimal and maximal spectra as the Stone-\v{C}ech compactification}

\author[A. Tarizadeh, M.R. Rezaee]{Abolfazl Tarizadeh, Mohammad Reza Rezaee}
\address{Department of Mathematics, Faculty of Basic Sciences, University of Maragheh \\
P. O. Box 55136-553, Maragheh, Iran.}
\email{ebulfez1978@gmail.com \\
atarizadeh@maragheh.ac.ir}

\address{Department of Mathematics, Faculty of Basic Sciences, University of Maragheh \\
P. O. Box 55136-553, Maragheh, Iran.}
\email{mohammadrezae.rezaee@yahoo.com}

\date{}
\subjclass[2010]{14A05, 54D35, 14M27, 54A20, 13A15}
\keywords{ Stone-\v{C}ech compactification; minimal spectrum; maximal spectrum; power set ring; Zariski convergent; ultra-ring; Alexandroff compactification \\}

\begin{abstract}  In this paper, new advances on the compactifications of topological spaces, especially on the Stone-\v{C}ech and Alexandroff compactifications have been made. Among the main results, it is proved that the minimal spectrum of the direct product of a family of integral domains indexed by a set $X$ is the Stone-\v{C}ech compactification of the discrete space $X$. Dually, it is proved that the maximal spectrum of the direct product of a family of local rings indexed by $X$ is also the Stone-\v{C}ech compactification of the discrete space $X$. The Alexandroff (one-point) compactification of a discrete space is constructed by a new method. Next, we proceed to give a natural and quite simple way to construct ultra-rings. Then this new approach is used to obtain several new results on the Stone-\v{C}ech compactification.
\end{abstract}

\maketitle

\section{Introduction}

Compactification is one of the main topics which is investigated in this paper from a purely algebraic perspective. Among various compactifications, the Stone-\v{C}ech compactification of a discrete space $X$ is particularly important. One of the main reasons of its importance is that it admits a semigroup structure whenever $X$ is a semigroup, and this semigroup structure has vast and interesting applications in diverse fields of mathematics specially in combinatorial number theory, Ramsey theory, topological dynamics and Ergodic theory. An accessible concrete description of this compactification often remains elusive. For instance the semigroup $\beta\mathbb{N}$, the Stone-\v{C}ech compactification of the natural numbers, is amazingly complicated and there are some unanswered questions about its semigroup structure. For example, whether or not $\beta\mathbb{N}$ contains any elements of finite order which are not idempotent still remains a challenging open problem. See \cite{Hindman-Strauss} and \cite{Walker} and their rich bibliography for further studies. Perhaps as another main reason for the importance of the Stone-\v{C}ech compactification of a discrete space is its vital role in proving Theorem \ref{Theorem niceTS} which asserts that every topological space admits the Stone-\v{C}ech compactification. \\

Classically, the Stone-\v{C}ech compactification of a discrete space is usually constructed via the ultrafilters of that space. In this paper, we find two new and interesting ways to construct this compactification using only the standard and elementary methods of commutative algebra. In fact in Theorem \ref{Theorem STIII}, we prove that the minimal spectrum of the direct product of a family of integral domains indexed by a set $X$ is the Stone-\v{C}ech compactification of the discrete space $X$. In Theorem \ref{Theorem STI}, it is shown that the maximal spectrum of the direct product of a family of local rings indexed by $X$ is the Stone-\v{C}ech compactification of the discrete space $X$. These results improve all of the former constructions of the Stone-\v{C}ech compactification of a discrete space, and also show that this compactification is independent of choosing of integral domains and local rings. In particular, we get that $\beta X=\Spec\mathcal{P}(X)$. The classical construction is also recovered (see Remark \ref{Remark II}). Throughout this paper, $\beta X$ denotes the Stone-\v{C}ech compactification of the discrete space $X$. These results allow us to understand the number of prime ideals of the infinite direct products of integral domains and local rings more precisely. As another application, the Stone-\v{C}ech compactification of an arbitrary topological space $X$ is deduced from the Stone-\v{C}ech compactification of the discrete space $X$ by passing to a certain quotient (see Theorem \ref{Theorem niceTS}). It is worth mentioning that our results considerably generalize several related results in the literature (see e.g. \cite{Aoyama}). \\

In \S6, using ultra-rings and Theorems \ref{Theorem STIII} and \ref{Theorem STI}, then we obtain new results on the Stone-\v{C}ech compactification (see Theorems \ref{Theorem Ultra 1} and \ref{Theorem Ultra 2}). \\

We introduce a new way to build the Alexandroff (one-point) compactification of a discrete space, see Corollary \ref{Corollary IV}. This result tells us that for any set $X$, then $\alpha X=\Spec(\mathcal{R})$. Here $\alpha X$ denotes the Alexandroff compactification of the discrete space $X$ and $\mathcal{R}$ is a certain subring of $\mathcal{P}(X)$. Then in Theorem \ref{Theorem IV}, we show that every totally disconnected compactification of a discrete space $X$ is precisely of the form $\Spec(\mathcal{R}')$ where the ring $\mathcal{R}'$ satisfies in the extensions of rings $\mathcal{R}\subseteq\mathcal{R}'\subseteq\mathcal{P}(X)$. After proving this result, we were informed that it is also proved in \cite[Theorems 2.2 and 2.3]{Swamy et al.} by another approach using Boolean algebras. In summary, our result shows that all of the totally disconnected compactifications of a discrete space are in the scope of the Zariski topology. In particular this class, up to isomorphisms, forms a set and the extensions of the corresponding rings put a partial order over this set in a way that the Alexandroff compactification is the minimal one and the Stone-\v{C}ech compactification is the maximal one. \\

It is well known that if the discrete space $X$ is also a (commutative) semigroup then its operation can be extended uniquely to an operation on $\beta X$ which forms a semigroup structure as well, see \cite[Theorems 4.1 and 4.4]{Hindman-Strauss}. This result opens new horizons to explore the basic and also sophisticated properties of the semigroup $\beta X$. Although some of them have been done in the literature over the years (see \cite{Hindman-Strauss} and its bibliography), there is a pressing need for new constructions to aid the development and the understanding the algebraic structure of this semigroup specially $\beta\mathbb{N}$ more deeply. We have made very little contributions to this subject but the results are sufficiently general (Theorems \ref{Theorem DV}, \ref{Theorem the4ta} and \ref{Proposition Elman-Asena}). Indeed, in Theorem \ref{Theorem DV}, we reformulate this important result into a more standard form and then it is proven by a new approach. Then in Theorems \ref{Theorem the4ta} and \ref{Proposition Elman-Asena}, various aspects of the semigroup $\beta X$ are investigated, specially it is shown that this semigroup structure is actually functorial. Finally, in Section 10, the absolutely flatness of the total ring of fractions is investigated.

\section{Preliminaries}

In this paper, all rings are commutative.  If $\phi:R\rightarrow R'$ is a morphism of rings then the induced map $\Spec(R')\rightarrow\Spec(R)$ given by $\mathfrak{p}\rightsquigarrow\phi^{-1}(\mathfrak{p})$ is
denoted by $\Spec(\phi)$, (sometimes it is also denoted by $\phi^{\ast}$). Let $X$ be a set and $I$ an ideal of the power set ring $\mathcal{P}(X)$. If $A,B\in I$ then $A\cup B=A+B+A\cap B\in I$. Also if $A\in I$ and $B\subseteq A$, then $B\in I$. For the definition of power set ring see e.g. \cite[\S2]{Tarizadeh}.
By $\Fin(X)$ we mean the set of all finite subsets of $X$, it is an ideal of $\mathcal{P}(X)$. By $\Clop(X)$ we mean the set of clopen (both open and closed) subsets of a topological space $X$ (for more information see \cite[\S3]{Tarizadeh}). By a compact space we mean a quasi-compact and Hausdorff topological space. \\

Let $R$ be a ring. The set of minimal primes of $R$ is denoted by $\Min(R)$ and the set of maximal ideals of $R$ is denoted by $\Max(R)$. Note that the induced Zariski topology over $\Min(R)$ is not necessarily quasi-compact. The Jacobson radical of $R$ is denoted by $\mathfrak{J}$.

\begin{definition} The Stone-\v{C}ech compactification of a topological space $X$ is the pair $(\beta X, \eta)$ where $\beta X$ is a compact space and $\eta:X\rightarrow\beta X$ is a continuous map such that the following universal property holds. For each such pair $(Y,\phi)$, i.e. $Y$ is a compact space and $\phi:X\rightarrow Y$ is a continuous map, then there exists a unique continuous map $\widetilde{\phi}:\beta X\rightarrow Y$ such that $\phi=\widetilde{\phi}\circ\eta$.
\end{definition}

All of the remaining undefined notions such as the flat topology, retraction, mp-ring and etc can be found in \cite{Tarizadeh-Aghajani}, \cite{A. Tarizadeh ftnm}, \cite{Abolfazl}, \cite{Tarizadeh 303}, \cite{A-Tarizadeh-Racsam} and \cite{A. Tarizadeh Racsam2}.

\section{Minimal spectrum as Stone-\v{C}ech compactification}

The main result of this section (Theorem \ref{Theorem STIII}) asserts that the minimal spectrum of the direct product of a family of integral domains indexed by a set $X$ is the Stone-\v{C}ech compactification of the discrete space $X$. \\

We start with the following result which generalizes \cite[p. 460]{Marco-Orsatti}.

\begin{proposition}\label{Proposition IIST41} Consider the canonical ring map $\pi:R\rightarrow S^{-1}R$ where $S$ is a multiplicative subset of a ring $R$, and let $f\in R$. Then $f\in\bigcap\limits_{\mathfrak{p}\in\Ima\pi^{\ast}}\mathfrak{p}$ if and only if there exists some $g\in S$ such that $fg$ is nilpotent.
\end{proposition}

\begin{proof} If $f\in\bigcap\limits_{\mathfrak{p}\in\Ima\pi^{\ast}}\mathfrak{p}$ then $f/1\in\bigcap\limits_{\mathfrak{p}\in\Ima\pi^{\ast}}S^{-1}\mathfrak{p}=
\bigcap\limits_{\mathfrak{q}
\in\Spec(S^{-1}R)}\mathfrak{q}=\sqrt{0}$. Thus there exist a natural number $n\geqslant1$ and some $g\in S$ such that $f^{n}g=0$. Hence, $fg$ is nilpotent. The reverse implication is easy.
\end{proof}

\begin{corollary}$($\cite[Lemma 1.1]{Henriksen-Jerison} and \cite[Lemma 3.1]{Kist}$)$\label{Corollary minimal primes} Let $\mathfrak{p}$ be a prime ideal of a ring $R$. Then $\mathfrak{p}$ is a minimal prime of $R$ if and only if for each $f\in\mathfrak{p}$ there exists some $g\in R\setminus\mathfrak{p}$ such that $fg$ is nilpotent.
\end{corollary}

\begin{proof} It is an immediate consequence of Proposition \ref{Proposition IIST41}.
\end{proof}

\begin{theorem}\label{Theorem 2000st} Let $R$ be a ring. Then the induced Zariski topology over $\Min(R)$ is finer than the induced flat topology. These two topologies over $\Min(R)$ are the same if and only if $\Min(R)$ is Zariski compact.
\end{theorem}

\begin{proof} Let $f\in R$. If $\mathfrak{p}\in W=\Min(R)\cap V(f)$ then by Corollary \ref{Corollary minimal primes}, there exists some $g\in R\setminus\mathfrak{p}$ such that $fg$ is nilpotent. This yields that $\mathfrak{p}\in\Min(R)\cap D(g)\subseteq W$. Therefore $W$ is a Zariski open of $\Min(R)$. Hence the Zariski topology over $\Min(R)$ is finer than the flat topology. For any ring $R$, then by \cite[Lemma 3.2]{Tarizadeh-Aghajani}, $\Min(R)$ is Zariski Hausdorff. Also $\Min(R)$ is flat quasi-compact. Therefore if these two topologies over $\Min(R)$ are the same then $\Min(R)$ is Zariski compact. Conversely, suppose $\Min(R)$ is Zariski compact. In the above we observed that $U=\Min(R)\cap D(f)$ is a Zariski clopen of $\Min(R)$. Every closed subspace of a quasi-compact space is quasi-compact. Thus there exist finitely many elements $g_{1},\ldots,g_{n}\in R$ such that $U^{c}=\Min(R)\setminus U=\bigcup\limits_{i=1}^{n}\Min(R)\cap D(f_{i})$. It follows that $U=\Min(R)\cap V(I)$ where $I=(g_{1},\ldots,g_{n})$ is a finitely generated ideal of $R$. Thus $U$ is a flat open of $\Min(R)$.
\end{proof}

Throughout this paper, $\Lambda=\prod\limits_{x\in X}R_{x}$ where each $R_{x}$ is an integral domain. For each $f=(f_{x})\in\Lambda$, the set $\Supp(f)=\{x\in X: f_{x}\neq0\}$ is simply denoted by $\Su(f)$. Clearly $\Su(fg)=\Su(f)\cap\Su(g)$ for all $f,g\in\Lambda$.

\begin{corollary}\label{Lemma STIII} The space $\Min(\Lambda)$ is Zariski compact.
\end{corollary}

\begin{proof} By Theorem \ref{Theorem 2000st}, it suffices to show that for each $f\in\Lambda$ then $U=\Min(\Lambda)\cap D(f)$ is a flat open of $\Min(\Lambda)$. Consider the sequence $e=(e_{x})\in\Lambda$ where $e_{x}$ is either $0$ or $1$, according as $x\in\Su(f)$ or $x\notin\Su(f)$. Then clearly $ef=0$ and $g=ge$ for all $g\in\Ann(f)$. Hence $\Ann(f)$ is generated by the sequence $e$. Now let $\mathfrak{p}\in\Min(\Lambda)\cap V(e)$. If $f\in\mathfrak{p}$ then by Corollary \ref{Corollary minimal primes}, there exists some $h\in\Lambda\setminus\mathfrak{p}$ such that $fh$ is nilpotent. But $\Lambda$ is a reduced ring. Hence $h\in\Ann(f)$. Thus $h=he\in\mathfrak{p}$. But this is a contradiction. This shows that $U=\Min(\Lambda)\cap V(e)$ is a flat open of $\Min(\Lambda)$.
\end{proof}

For each $x\in X$ then $\mathfrak{p}_{x}:=\Ker\pi_{x}$ is a minimal prime of $\Lambda$ and it is generated by the sequence $1-\Delta_{x}$ where $\pi_{x}:\Lambda\rightarrow R_{x}$ is the canonical projection, $\Delta_{x}=(\delta_{x,y})_{y\in X}$ and $\delta_{x,y}$ is the Kronecker delta. \\

Now we are ready to prove the main result of this section:

\begin{theorem}\label{Theorem STIII} The space $\Min(\Lambda)$ together with the canonical map $\eta:X\rightarrow\Min(\Lambda)$ given by $x\rightsquigarrow\mathfrak{p}_{x}$ is the Stone-\v{C}ech compactification of the discrete space $X$.
\end{theorem}

\begin{proof} By Corollary \ref{Lemma STIII}, the space $\Min(\Lambda)$ is compact. It remains to check out the universal property of the Stone-\v{C}ech compactification. Let $Y$ be a compact topological space and $\phi:X\rightarrow Y$ a function. We shall find a continuous function $\widetilde{\phi}:\Min(\Lambda)\rightarrow Y$ such that $\phi=\widetilde{\phi}\circ\eta$ and then we show that such function is unique. If $\mathfrak{p}\in\Min(\Lambda)$ then the subsets $\Su(f)$ with $f\in\Lambda\setminus\mathfrak{p}$ have the finite intersection property. It follows that the subsets $\phi\big(\Su(f)\big)$ and so their closures $\overline{\phi\big(\Su(f)\big)}$ with $f\in\Lambda\setminus\mathfrak{p}$ have the finite intersection property. This yields that $\bigcap\limits_{f\in \Lambda\setminus\mathfrak{p}}\overline{\phi\big(\Su(f)\big)}\neq\emptyset$ because $Y$ is quasi-compact. We claim that this intersection has exactly one point. If $y$ and $y'$ are two distinct points of the intersection then there exist disjoint opens $U$ and $V$ in $Y$ such that $y\in U$ and $y'\in V$. Then consider the sequence $f\in\Lambda$ where $f_{x}$ is either $0$ or $1$, according as $x\in\phi^{-1}(U)$ or $x\notin\phi^{-1}(U)$. Then we have either $f\in\mathfrak{p}$ or $1-f\in\mathfrak{p}$ since $f$ is an idempotent. If $f\in\mathfrak{p}$ then $\phi^{-1}(V)\cap\Su(1-f)\neq\emptyset$. So we may choose some $x$ in this intersection. Thus $x\notin\phi^{-1}(U)$, hence $f_{x}=1$. But this is a contradiction since $x\in\Su(1-f)$. If $1-f\in\mathfrak{p}$ then $\phi^{-1}(U)\cap\Su(f)\neq\emptyset$, but this is again a contradiction. Hence, there exists a unique point $y_{\mathfrak{p}}\in Y$ such that $\bigcap\limits_{f\in \Lambda\setminus\mathfrak{p}}\overline{\phi\big(\Su(f)\big)}=
\{y_{\mathfrak{p}}\}$. This establishes the claim.
Then we define the map $\widetilde{\phi}:\Min(\Lambda)\rightarrow Y$ as $\mathfrak{p}\rightsquigarrow y_{\mathfrak{p}}$. It is easy to see that $\phi(x)\in\bigcap\limits_{f\in \Lambda\setminus\mathfrak{p}_{x}}\phi\big(\Su(f)\big)$ for all $x\in X$. Therefore $\phi=\widetilde{\phi}\circ\eta$. Now we show that $\widetilde{\phi}$ is continuous. Let $U$ be an open of $Y$ and let $\mathfrak{p}\in(\widetilde{\phi})^{-1}(U)$. There exists an open neighborhood $V$ of $y_{\mathfrak{p}}$ such that $\overline{V}\subseteq U$, because it is well known that every compact space is a normal space. Let $h\in \Lambda$ be a sequence which is defined as $h_{x}=1$ or $h_{x}=0$, according as $x\in\phi^{-1}(V)$ or $x\notin\phi^{-1}(V)$. Then $\mathfrak{p}\in\D(h)$, since if $h\in\mathfrak{p}$ then $1-h\notin\mathfrak{p}$ and so $\phi^{-1}(V)\cap\Su(1-h)\neq\emptyset$, which is impossible. To conclude the continuity of $\widetilde{\phi}$ we show that $\Min(\Lambda)\cap\D(h)\subseteq(\widetilde{\phi})^{-1}(U)$. Suppose there exists some $\mathfrak{q}\in\Min(\Lambda)\cap\D(h)$ such that $y_{\mathfrak{q}}\notin U$. Thus $y_{\mathfrak{q}}\in W:=Y\setminus\overline{V}$. It follows that $W\cap\phi\big(\Su(h)\big)\neq\emptyset$. But this is impossible since $\Su(h)=\phi^{-1}(V)$ and so $W\cap\phi\big(\Su(h)\big)\subseteq W\cap V=\emptyset$. Therefore $\widetilde{\phi}$ is continuous. If $\Min(\Lambda)\cap D(f)$ is non-empty then $f\neq0$ and so there exists some $x\in X$ such that $\mathfrak{p}_{x}\in D(f)$. This shows that
$\eta(X)$ is a dense subspace of $\Min(\Lambda)$, hence the uniqueness of $\widetilde{\phi}$ is deduced from the basic fact that if two continuous maps into a Hausdorff space agree on a dense subspace of the domain, they are equal.
\end{proof}

\begin{lemma}\label{Lemma DII} If each $R_{x}$ is a field, then every prime ideal of $\Lambda$ is a maximal ideal.
\end{lemma}

\begin{proof} Let $\mathfrak{p}$ be a prime ideal of $\Lambda$ and $f\in \Lambda\setminus\mathfrak{p}$. Then consider the sequence $g=(g_{x})\in\Lambda$ where $g_{x}$ is $1$ or $1/f(x)$, according as $f_{x}=0$ or $f_{x}\neq0$. Then it is obvious that $f(1-fg)=0\in\mathfrak{p}$. This yields that $1-fg\in\mathfrak{p}$. Therefore $\Lambda/\mathfrak{p}$ is a field. As a second proof, the assertion is also deduced from the fact that $\Lambda$ is an absolutely flat ring.
\end{proof}

\begin{corollary}\label{Theorem DI} The space $\Spec(\Lambda)$ together with the canonical map $\eta:X\rightarrow\Spec(\Lambda)$
is the Stone-\v{C}ech compactification of the discrete space $X$ if and only if each $R_{x}$ is a field.
\end{corollary}

\begin{proof} If each $R_{x}$ is a field then the assertion is deduced from Theorem \ref{Theorem STIII} and Lemma \ref{Lemma DII}. Conversely, if $\mathfrak{m}$ is a maximal ideal of $R_{x}$ then $\pi^{-1}_{x}(\mathfrak{m})=\pi^{-1}_{x}(0)$ because $\Spec(\Lambda)$ is Hausdorff and so every prime ideal of $\Lambda$ is a maximal ideal. But $\pi_{x}$ is surjective and so the induced map $\pi^{\ast}_{x}$ is injective. Therefore the zero ideal of $R_{x}$ is a maximal ideal and so it is a field.
\end{proof}

\begin{corollary}\label{Corollary STII} The space $\Spec\mathcal{P}(X)$ together with the canonical map $\eta:X\rightarrow\Spec\mathcal{P}(X)$ given by $x\rightsquigarrow\mathfrak{m}_{x}=\mathcal{P}(X\setminus\{x\})$
is the Stone-\v{C}ech compactification of the discrete space $X$.
\end{corollary}

\begin{proof} The map $\mathcal{P}(X)\rightarrow\prod\limits_{x\in X}\mathbb{Z}_{2}$ given by $A\rightsquigarrow\chi_{A}$ is an isomorphism of rings where $\chi_{A}$ is the characteristic function of $A$ and $\mathbb{Z}_{2}=\{0,1\}$. Then apply Corollary \ref{Theorem DI}.
\end{proof}

\begin{remark}\label{Remark II} Here we establish a bridge that allows us to translate all of the theory of Boolean algebras into the standard language of commutative algebra (and vice versa). For instance, the classical approach to construct the Stone-\v{C}ech compactification of a discrete space $X$ is easily recovered. Indeed, if $X$ is a set then one can easily check that the map $M\rightsquigarrow \mathcal{P}(X)\setminus M=\{A\in\mathcal{P}(X): A^{c}\in M\}$ is a homeomorphism from $\Spec\mathcal{P}(X)$ onto $\mathscr{F}(X)$, the space of ultrafilters on $X$ equipped with the Stone topology. Recall that the collection of $d(A)=\{F\in\mathscr{F}(X):A\in F\}$ with $A\in\mathcal{P}(X)$ forms a base for the opens of the Stone topology. The space $\mathscr{F}(X)$
is called the Stone space of the Boolean algebra $\mathcal{P}(X)$. Note that the above identification can be generalized to any Boolean ring $R$. In fact, the map $M\rightsquigarrow R\setminus M$ is a homeomorphism from $\Spec(R)$ onto the Stone space of the corresponding Boolean algebra of $R$.
\end{remark}

Remember that for any two objects $X$ and $Y$ of a category $\mathscr{C}$, by $\Mor_{\mathscr{C}}(X,Y)$ we mean the set of all morphisms of $\mathscr{C}$ from $X$ to $Y$.

\begin{corollary}\label{Corollary DI} For any two sets $X$ and $Y$ then we have the following canonical bijections: $$\Mor_{\Set}(X,\beta Y)\simeq\Mor_{\Top}(\beta X, \beta Y)\simeq\Mor_{\Ring}\big(\mathcal{P}(Y),\mathcal{P}(X)\big).$$
\end{corollary}

\begin{proof} The first bijection follows form Corollary \ref{Corollary STII}, and the second bijection is an immediate consequence of \cite[Theorem 5.6]{Tarizadeh}.
\end{proof}

\section{The Stone-\v{C}ech compactification of an arbitrary space}

In this section we give a new proof to the fact that every topological space admits the Stone-\v{C}ech compactification. To realize this goal, we first obtain some results which are interesting in their own right. We should mention that these results are not new and can be found in the literature which are obtained by using the theory of ultrafilters. Maybe the only novelty is that we will use only the ring-theoretical methods, it seems that this approach is much simpler than the theory of ultrafilters. \\

We begin with the following key definition (this notion is due to Henri Cartan, we interpreted it into the language of commutative algebra).

\begin{definition} Let $X$ be a topological space, $x\in X$ and $M$ a maximal ideal of $\mathcal{P}(X)$. We say that $M$ is convergent to the point $x$ if whenever $U$ is an open subset of $X$ containing $x$, then $M\in D(U)$.
\end{definition}

\begin{lemma}\label{Lemma TTPN0923465} Let $X$ be a set. If $M$ is a maximal ideal of $\mathcal{P}(X)$ then $\mathcal{P}(\eta)^{\ast}(M)$ is convergent to the point $M\in \beta X=\Spec\mathcal{P}(X)$.
\end{lemma}

\begin{proof} Let $U$ be an open of $\beta X$ such that $M\in U$. If $U\in\mathcal{P}(\eta)^{\ast}(M)$ then $\eta^{-1}(U)\in M$. But there exists some $A\in\mathcal{P}(X)$ such that $M\in D(A)\subseteq U$. If $x\in A$ then $\eta(x)=\mathfrak{m}_{x}\in D(A)$ and so $x\in\eta^{-1}(U)$. This shows that $A\subseteq\eta^{-1}(U)$. Thus $A\in M$. But this is a contradiction.
\end{proof}

\begin{lemma}\label{Lemma BPNT} Let $X$ be a topological space and let $A$ be a subset of $X$ with the property that $D(A)$ contains every maximal ideal of $\mathcal{P}(X)$ which is convergent to a point of $A$. Then $A$ is an open subset of $X$.
\end{lemma}

\begin{proof} Take $x\in A$ and let $\mathcal{S}$ be the set of all opens of $X$ which are containing $x$. Then by the hypothesis, the ideal of $\mathcal{P}(X)$ generated by $A$ and the elements $U^{c}=X\setminus U$ with $U\in\mathcal{S}$ is the whole ring. Thus we may find a finite number $U_{1},...,U_{n}$ of elements of $\mathcal{S}$ such that $X=A\cup(\bigcup\limits_{i=1}^{n}U^{c}_{i})$. It follows that $x\in\bigcap\limits_{i=1}^{n}U_{i}\subseteq A$. Hence, $A$ is an open of $X$.
\end{proof}

Note that the converse of the above lemma holds trivially. \\

Let $\phi:X\rightarrow Y$ be a continuous map of topological spaces. If a maximal ideal $M$ of $\mathcal{P}(X)$ converges to some point $x\in X$, then clearly $\mathcal{P}(\phi)^{\ast}(M)$ is convergent to $\phi(x)$. In the following result we establish its converse.

\begin{corollary}\label{Lemma TPosD} Let $\phi:X\rightarrow Y$ be a function between topological spaces with the property that $\mathcal{P(\phi)^{\ast}}(M)$ is convergent to $\phi(x)$ whenever a maximal ideal $M$ of $\mathcal{P}(X)$ converges to some point $x\in X$. Then $\phi$ is continuous.
\end{corollary}

\begin{proof} It is easily deduced from Lemma \ref{Lemma BPNT}.
\end{proof}

Now we establish the main result of this section. This result shows that the Stone-\v{C}ech construction can be performed for any topological space $X$, but in that case the canonical map from $X$ to its Stone-\v{C}ech compactification need not be a homeomorphism onto its image (and sometimes is not even injective).

\begin{theorem}\label{Theorem niceTS} Every topological space $X$ admits the Stone-\v{C}ech compactification.
\end{theorem}

\begin{proof} Consider the equivalence relation $\sim$ on $\beta X=\Spec\mathcal{P}(X)$ defined as $M\sim N$ if $\phi:X\rightarrow Y$ is a continuous function to a compact space $Y$ then $\widetilde{\phi}(M)=\widetilde{\phi}(N)$ where $\widetilde{\phi}:\beta X\rightarrow Y$ is the unique continuous function such that $\phi=\widetilde{\phi}\circ\eta$, see the proof of Theorem \ref{Theorem STIII}.  Now to prove that the pair $(X', \pi\circ\eta)$ is the Stone-\v{C}ech compactification of the space $X$ it suffices to show that $\pi\circ\eta:X\rightarrow X'$ is continuous where $\pi:\beta X\rightarrow X'=\beta X/\sim$ is the canonical map and $X'$ is equipped with the quotient topology. To prove the continuity of $\pi\circ\eta$, by Corollary \ref{Lemma TPosD}, it will be enough to show that if a maximal ideal $M$ of $\mathcal{P}(X)$ converges to some point $x\in X$ then $\mathcal{P}(\pi\circ\eta)^{\ast}(M)$ is convergent to the point $(\pi\circ\eta)(x)$. We have  $\mathcal{P}(\pi\circ\eta)^{\ast}(M)=
\mathcal{P}(\pi)^{\ast}\big(\mathcal{P}(\eta)^{\ast}(M)\big)$. By Lemma \ref{Lemma TTPN0923465}, $N:=\mathcal{P}(\eta)^{\ast}(M)$ is convergent to the point $M\in \beta X$. Thus $\mathcal{P}(\pi)^{\ast}(N)$ is convergent to the point $\pi(M)$ since $\pi$ is continuous. Then we show that $M\sim\mathfrak{m}_{x}$. Because take $A\in\mathcal{P}(X)\setminus M$ and let $V$ be an open of a compact space $Y$ such that $\phi(x)\in V$ where $\phi:X\rightarrow Y$ is a continuous map. Then $\phi^{-1}(V)\notin M$. Note that $\Su(A)=A$. Now if $V\cap \phi(A)=\emptyset$ then $A\in M$, a contradiction. Hence, $\phi(x)\in\overline{\phi\big(\Su(A)\big)}$. Thus by the definition of $\widetilde{\phi}$, see the proof of Theorem \ref{Theorem STIII}, we get that $\phi(x)=\widetilde{\phi}(M)$
and so $M\sim\mathfrak{m}_{x}$. Therefore $\mathcal{P}(\pi\circ\eta)^{\ast}(M)$ is convergent to the point $\pi(M)=(\pi\circ\eta)(x)$. Note that during to verify the universal property of the Stone-\v{C}ech compactification for the pair $(X', \pi\circ\eta)$, the uniqueness is deduced from the fact that $(\pi\circ\eta)(X)$ is a dense subspace of $X'$.
\end{proof}

Recall that by a \emph{compactification} of a topological space $X$ we mean a compact space $\widetilde{X}$ together with an open embedding (a continuous injective open map) $\eta:X\rightarrow\widetilde{X}$ such that $\eta(X)$ is a dense subspace of $\widetilde{X}$. If moreover, $\widetilde{X}\setminus\eta(X)$ consisting only a single point then $\widetilde{X}$ is called the \emph{one-point} or the \emph{Alexandroff compactification of $X$} and it is often denoted by $\alpha X$, and this single point is called the \emph{point at infinity of $X$}. \\

\begin{remark} Note that the map $\eta:X\rightarrow\Min(\Lambda)$ given by $x\rightsquigarrow \mathfrak{p}_{x}$ is an open embedding, because $\{\mathfrak{p}_{x}\}=\Min(\Lambda)\cap D(\Delta_{x})$ for all $x\in X$. Hence, $\Min(\Lambda)$ is also a compactification of the discrete space $X$ in the above sense. But it is important to notice that the Stone-\v{C}ech compactification of an arbitrary topological space need not be a compactification in the above sense. In fact, it is well known that the canonical map from a topological space $X$ to its Stone-\v{C}ech compactification induces a homeomorphism onto its image if and only if $X$ is a Tychonoff space. Thus for general space $X$, this map need not be injective. It is also well known that the canonical map from a topological space $X$ to its Stone-\v{C}ech compactification is an open embedding if and only if $X$ is locally compact.
\end{remark}

\section{Maximal spectrum as Stone-\v{C}ech compactification}

The main result of this section (Theorem \ref{Theorem STI}) asserts that the maximal spectrum of the direct product of a family of local rings indexed by a set $X$ is the Stone-\v{C}ech compactification of the discrete space $X$. Then some applications are also given. \\

Let $R$ be a ring and $f\in R$. If $\mathfrak{m}\in U=\Max(R)\cap D(f)$ then there exist some $g\in\mathfrak{m}$ and $h\in R$ such that $1=fh+g$. This yields that $\mathfrak{m}\in\Max(R)\cap V(g)\subseteq U$. Thus $U$ is a flat open of $\Max(R)$. Therefore the induced flat topology over $\Max(R)$ is finer than the induced Zariski topology.

\begin{proposition}\label{Proposition IST} For a ring $R$ the following statements are equivalent. \\
$\mathbf{(i)}$ $R/\mathfrak{J}$ is a zero dimensional ring. \\
$\mathbf{(ii)}$ The induced Zariski and flat topologies over $\Max(R)$ are the same. \\
$\mathbf{(iii)}$ $\Max(R)$ is flat compact.
\end{proposition}

\begin{proof} $\mathbf{(i)}\Rightarrow\mathbf{(ii)}:$ If $f\in R$ then there exists some $g\in R$ such that $f(1-fg)\in\mathfrak{J}$, because $R/\mathfrak{J}$ is reduced and so it is absolutely flat. It follows that $\Max(R)\cap V(f)=\Max(R)\cap D(1-fg)$. \\
$\mathbf{(ii)}\Rightarrow\mathbf{(iii)}:$ The subset $\Max(R)$ is Zariski quasi-compact and flat Hausdorff. \\
$\mathbf{(iii)}\Rightarrow\mathbf{(i)}:$ See \cite[Theorem 4.5]{Tarizadeh 303}.
\end{proof}

\begin{lemma}\label{Lemma STII} Let $R$ be a ring such that $R/\mathfrak{J}$ is a zero dimensional ring. Then the clopens of $\Max(R)$ are precisely of the form $\Max(R)\cap V(f)$ where $f\in R$.
\end{lemma}

\begin{proof} By Proposition \ref{Proposition IST}, the Zariski and flat topologies over $\Max(R)$ are the same. If $f\in R$ then we observed that $\Max(R)\cap V(f)$ is a clopen of $\Max(R)$. Conversely, let $U$ be a clopen of $\Max(R)$. It is easy to see that every closed subspace of a quasi-compact space is quasi-compact. Hence, we may write $U=\bigcup\limits_{k=1}^{n}\Max(R)\cap V(I_{k})$ where each $I_{k}$ is a (finitely generated) ideal of $R$. This yields that $U=\Max(R)\cap V(I)$ where $I=I_{1}...I_{n}$. Similarly we get that $U^{c}=\Max(R)\setminus U=\Max(R)\cap V(J)$ where $J$ is a (finitely generated) ideal of $R$. It follows that $I+J=R$. Thus there exist some $f\in I$ and $g\in J$ such that $f+g=1$. This implies that $U=\Max(R)\cap V(f)$.
\end{proof}

Throughout this paper, $\Gamma=\prod\limits_{x\in X}R_{x}$ where each $R_{x}$ is a local ring with the maximal ideal $\mathfrak{m}_{x}$. For each $x\in X$ then $\mathfrak{M}_{x}:=\pi^{-1}_{x}(\mathfrak{m}_{x})$ is a maximal ideal of $\Gamma$, because the ring map $\Gamma/\mathfrak{M}_{x}\rightarrow R_{x}/\mathfrak{m}_{x}$ induced by the canonical projection $\pi_{x}:\Gamma\rightarrow R_{x}$ is an isomorphism. If $f=(f_{x})\in\Gamma$ then we define $\Omega(f)=\{x\in X: f_{x}\notin\mathfrak{m}_{x}\}$. It is obvious that $f$ is invertible in $\Gamma$ if and only if $\Omega(f)=X$. It is also easy to see that $\Omega(fg)=\Omega(f)\cap\Omega(g)$ for all $f,g\in\Gamma$.

\begin{lemma}\label{Lemma STI} Let $f\in\Gamma$. Then $\Omega(f)=\emptyset$ if and only if $f\in\mathfrak{J}$.
\end{lemma}

\begin{proof} If $\Omega(f)=\emptyset$ then $f_{x}\in\mathfrak{m}_{x}$ for all $x$. This yields that $\Omega(1+fg)=X$ for all $g\in\Gamma$. Thus $f\in\mathfrak{J}$. Conversely, if $f\in\mathfrak{J}$ then $f\in\mathfrak{M}_{x}$ for all $x$. So $\Omega(f)$ is empty.
\end{proof}

\begin{theorem}\label{Theorem STI} The space $\Max(\Gamma)$ together with the canonical map $\eta:X\rightarrow\Max(\Gamma)$ given by $x\rightsquigarrow\mathfrak{M}_{x}$ is the Stone-\v{C}ech compactification of the discrete space $X$.
\end{theorem}

\begin{proof} If $f\in\Gamma$ then consider the sequence $g=(g_{x})\in\Gamma$ such that $g_{x}$ is either $0$ or  $f^{-1}_{x}$, according as $f_{x}\in\mathfrak{m}_{x}$ or $f_{x}\notin\mathfrak{m}_{x}$. Then $\Omega(1+fh(1-fg))=X$ for all $h\in\Gamma$. Hence, $f(1-fg)\in\mathfrak{J}$.
Thus $\Gamma/\mathfrak{J}$ is absolutely flat. Therefore by Proposition \ref{Proposition IST}, the space $\Max(\Gamma)$ is compact. Then we verify the universal property of the Stone-\v{C}ech compactification. Let $Y$ be a compact topological space and $\phi:X\rightarrow Y$ a function. If $M\in\Max(\Gamma)$ then by Lemma \ref{Lemma STI}, the subsets $\Omega(f)$ with $f\in\Gamma\setminus M$ have the finite intersection property. Thus by a similar argument as applied in the proof of Theorem \ref{Theorem STIII}, there exists a unique point $y_{M}\in Y$ such that $\bigcap\limits_{f\in\Gamma\setminus M}\overline{\phi\big(\Omega(f)\big)}=
\{y_{M}\}$. Then we define the map $\widetilde{\phi}:\Max(\Gamma)\rightarrow Y$ as $M\rightsquigarrow y_{M}$. Again exactly like the proof of Theorem \ref{Theorem STIII}, it is shown that $\phi=\widetilde{\phi}\circ\eta$ and $\widetilde{\phi}$ is continuous. Finally, to prove the uniqueness of $\widetilde{\phi}$ it suffices to show that $\eta(X)$ is a dense subspace of $\Max(\Gamma)$. The space $\Max(\Gamma)$ is totally disconnected, see \cite[Proposition 4.4]{Tarizadeh 303}. It is well known that in a compact totally disconnected space, the collection of clopens is a base for the opens. Using this and Lemma \ref{Lemma STII}, then the collection of $\Max(\Gamma)\cap V(f)$ with $f\in\Gamma$ forms a base for the opens of $\Max(\Gamma)$. Now if $\Max(\Gamma)\cap V(f)$ is non-empty then $\Omega(f)\neq X$. Hence there exists some $x\in X$ such that $\mathfrak{M}_{x}\in\Max(\Gamma)\cap V(f)$. Therefore $\eta(X)$ is a dense subspace of $\Max(\Gamma)$.
\end{proof}

\begin{example} If $p$ is a prime number then for each $n\geqslant1$,  $\mathbb{Z}/p^{n}\mathbb{Z}$ is a local zero dimensional ring (its prime spectrum is a singleton), but the direct product ring $\prod\limits_{n\geqslant1}\mathbb{Z}/p^{n}\mathbb{Z}$ has infinite Krull dimension. This ring also has a huge number of prime ideals. In fact by Theorem \ref{Theorem STI}, the cardinality of its maximal ideals equals $2^{\mathfrak{c}}$ where $\mathfrak{c}$ is the cardinality of the continuum.
\end{example}

\begin{remark} The canonical map $\eta:X\rightarrow\Max(\Gamma)$ given by $x\rightsquigarrow\mathfrak{M}_{x}$ is an open embedding. In fact, $\{\mathfrak{M}_{x}\}=\Max(\Gamma)\cap D(\Delta_{x})$ for all $x\in X$. To see this let $M\in\Max(\Gamma)\cap D(\Delta_{x})$ and $f\in M$. If $f\notin\mathfrak{M}_{x}$ then $\Omega(1-\Delta_{x}+\Delta_{x}f)=X$ and so
$1-\Delta_{x}+\Delta_{x}f$ is invertible in the ring $\Gamma$. But this is a contradiction because $1-\Delta_{x}+\Delta_{x}f\in M$. Therefore $M\subseteq\mathfrak{M}_{x}$ and so $M=\mathfrak{M}_{x}$.
\end{remark}

\begin{corollary}\label{Corollary STI} There exists a unique homeomorphism: $$\xymatrix{\Min(\Lambda)\ar[r]^{\simeq\:\:}&\Max(\Gamma)}$$ such that $\mathfrak{p}_{x}$ is mapped into $\mathfrak{M}_{x}$ for all $x\in X$.
\end{corollary}

\begin{proof} It is deduced from the universal property of the Stone-\v{C}ech compactification by taking into account Theorems \ref{Theorem STIII} and \ref{Theorem STI}.
\end{proof}

In the next section, we will precisely determine the rule of isomorphism of Corollary \ref{Corollary STI}.

\begin{corollary} Let $R$ be a ring and let $X$ be a subset of $\Spec(R)$. Then the following spaces are canonically isomorphic $($up to a unique isomorphism$)$. \\
$\mathbf{(i)}$ $\Min(\prod\limits_{\mathfrak{p}\in X}R/\mathfrak{p})$. \\
$\mathbf{(ii)}$ $\Spec\big(\prod\limits_{\mathfrak{p}\in X}\kappa(\mathfrak{p})\big)$. \\
$\mathbf{(iii)}$ $\Max(\prod\limits_{\mathfrak{p}\in X}R_{\mathfrak{p}})$.
\end{corollary}

\begin{proof} It is an immediate consequence of Corollary \ref{Corollary STI}.
\end{proof}

If $X$ is a set with the cardinality $\kappa$ and $\widetilde{X}$ is a compactification of the discrete space $X$, then by \cite[Tag 0909]{Johan} and assuming the generalized continuum hypothesis, we have $|\widetilde{X}|\in\{\kappa, 2^{\kappa},2^{2^{\kappa}}\}$.

\begin{corollary}\label{Corollary STIII} If $X$ is an infinite set with the cardinality $\kappa$, then $|\Min(\Lambda)|=|\Max(\Gamma)|=|\Spec\mathcal{P}(X)|=
2^{2^{\kappa}}$.
\end{corollary}

\begin{proof} It follows from Corollaries \ref{Corollary STII} and \ref{Corollary STI} and the fact that the cardinality of the Stone-\v{C}ech compactification of the infinite discrete space $X$ is equal to $2^{2^{\kappa}}$. To see the proof of this fact please consider \cite[Theorem 3.58]{Hindman-Strauss} or \cite[Theorem on page 71]{Walker}.
\end{proof}

\begin{corollary}\label{Corollary cardinal2} Let $X$ and $Y$ be two sets with the cardinalities $\kappa$ and $\lambda$, respectively. Then the number of all ring maps $\mathcal{P}(X)\rightarrow\mathcal{P}(Y)$ is either $\kappa^{\lambda}$ or $2^{\lambda2^{\kappa}}$, according as $X$ is finite or infinite.
\end{corollary}

\begin{proof} It is deduced from Corollaries \ref{Corollary DI} and \ref{Corollary STIII}.
\end{proof}

If $\kappa$ is an infinite cardinal, then $\lambda2^{\kappa}=\max\{\lambda,2^{\kappa}\}$. To see its proof apply Cantor's theorem and \cite[p. 162, Lemma 6R]{Enderton} which states that $\kappa\kappa=\kappa$.

\begin{corollary} Let $X$ be a set with the cardinality $\kappa$. Then the number of all ring maps $\mathcal{P}(X)\rightarrow\mathcal{P}(X)$ is either $\kappa^{\kappa}$ or $2^{2^{\kappa}}$, according as $\kappa$ is finite or infinite.
\end{corollary}

\begin{proof} It is an immediate consequence of Corollary \ref{Corollary cardinal2}.
\end{proof}

\section{Ultra-rings and their applications in compactification}

In this section, we introduce a new way to construct the ultraproduct of rings which considerably simplifies the existence method in the literature (see e.g. \cite{Becker et al.}, \cite{Denef-Lipshitz}, \cite{Erman et al.} and \cite{Schoutens}). Then we use this new approach to determine precisely the isomorphisms whose rules are already obtained in an implicit way (see e.g. Corollary \ref{Corollary STI}). \\

Let $(R_{x})$ be a family of rings indexed by a set $X$ and let $R=\prod\limits_{x\in X}R_{x}$ be their direct product ring. Let $M$ be a maximal ideal of $\mathcal{P}(X)$. Then it can be easily seen that $M^{\ast}=\{f\in R: \Su(f)\in M\}$ is an ideal of $R$, because clearly $\Su(0)=\emptyset\in M$ and so $0\in M^{\ast}$, also $\Su(f+g)\subseteq\Su(f)\cup\Su(g)$ and $\Su(fg)\subseteq\Su(f)\cap\Su(g)$ for all $f,g\in R$. We call the quotient ring $R/M^{\ast}$ the ultraproduct (or, ultra-ring) of the family $(R_{x})$ with respect to $M$. \\

It is interesting to notice the map $\phi:R\rightarrow\mathcal{P}(X)$ given by $f\rightsquigarrow\Su(f)$ is not a morphism of rings, since it is not additive, in fact $\Su(f)+\Su(g)\subseteq\Su(f+g)$. In spite of this, the inverse image of each ideal of $\mathcal{P}(X)$ under $\phi$ is an ideal of $R$. In particular, $M^{\ast}=\phi^{-1}(M)$. \\

If each $R_{x}$ is a non-zero ring, then $M^{\ast}$ is a proper ideal of $R$. We also have the following result which shows that some properties (statements in first-order logic) of rings are preserved under the formation of ultraproducts.

\begin{theorem}\label{Theorem IUltra} For a ring $R=\prod\limits_{x\in X}R_{x}$ the following assertions hold. \\
$\mathbf{(i)}$ If each $R_{x}$ is a field, then $R/M^{\ast}$ is a field. \\
$\mathbf{(ii)}$ If each $R_{x}$ is an integral domain, then $R/M^{\ast}$ is an integral domain. \\
$\mathbf{(iii)}$ If each $R_{x}$ is a local ring, then $R/M^{\ast}$ is a local ring. \\
$\mathbf{(iv)}$ If each $K_{x}$ is the fraction field of an integral domain $R_{x}$, then the ultra-ring of the family $(K_{x})$ with respect to $M$ is  the fraction field of $R/M^{\ast}$. \\
$\mathbf{(v)}$ If each $R_{x}$ is a local ring with the residue field $K_{x}$, then the ultra-ring of the family $(K_{x})$ with respect to $M$ is  the residue field of $R/M^{\ast}$.
\end{theorem}

\begin{proof} $\mathbf{(i)}:$ Take $f\in R\setminus M^{\ast}$ and consider the sequence $g=(g_{x})\in R$ where each $g_{x}$ is either $f^{-1}_{x}$ or $1$, according as $x\in\Su(f)$ or $x\notin\Su(f)$. Then clearly $\Su(1-fg)\subseteq\Su(f)^{c}\in M$. Thus $\Su(1-fg)\in M$ and so $1-fg\in M^{\ast}$. \\
$\mathbf{(ii)}:$ Suppose $fg\in M^{\ast}$ for some $f,g\in R$. Then clearly $\Su(f)\cap\Su(g)\subseteq\Su(fg)\in M$. Thus $\Su(f)\cap\Su(g)\in M$. It follows that either $f\in M^{\ast}$ or $g\in M^{\ast}$. \\
$\mathbf{(iii)}:$ Clearly $M^{\flat}=\{f\in R:\Omega(f)\in M\}$ is a proper ideal of $R$ and $M^{\ast}\subseteq M^{\flat}$, since $\Omega(f)=\{x\in X: f_{x}\notin\mathfrak{m}_{x}\}\subseteq\Su(f)$ for all $f\in R$ where  $\mathfrak{m}_{x}$ is the maximal ideal of $R_{x}$. If $f\in R\setminus M^{\flat}$ then $\Su(1-fg)\subseteq\Omega(f)^{c}\in M$ where $g=(g_{x})$ and each $g_{x}$ is either $f^{-1}_{x}$ or $1$, according as $x\in\Omega(f)$ or $x\notin\Omega(f)$. Therefore  $\Su(1-fg)\in M$ and so $1-fg\in M^{\ast}$. Hence, $M^{\flat}/M^{\ast}$ is the only maximal ideal of $R/M^{\ast}$. The proof of  $\mathbf{(iv)}$ is easy and left as an exercise to the reader. \\
$\mathbf{(v)}:$ It suffices to show that the map $R/M^{\flat}\rightarrow R'/M^{\ast}$ given by $f+M^{\flat}\rightsquigarrow \overline{f}+M^{\ast}$ is an isomorphism of rings where $R'=\prod\limits_{x\in X}K_{x}$, $M^{\ast}=\{g\in R': \Su(g)\in M\}$ and $\overline{f}=(f_{x}+\mathfrak{m}_{x})$ with $\mathfrak{m}_{x}$ is the maximal ideal of $R_{x}$. If $f\in M^{\flat}$ then $\Su(\overline{f})\subseteq\Omega(f)\in M$ and so $\Su(\overline{f})\in M$. Hence, the above map is well-defined. Clearly it is also an isomorphism.
\end{proof}

\begin{theorem}\label{Theorem Ultra 1} The map $\phi:\Spec\mathcal{P}(X)\rightarrow\Min(\Lambda)$ given by $M\rightsquigarrow M^{\ast}$ is a homeomorphism.
\end{theorem}

\begin{proof} First we need to show that $M^{\ast}$ is a minimal prime of $\Lambda$. By Theorem \ref{Theorem IUltra} (ii), $M^{\ast}$ is a prime ideal of $\Lambda$. Suppose there exists a prime ideal $\mathfrak{p}$ of $\Lambda$ such that $\mathfrak{p}\subseteq M^{\ast}$. If $f\in M^{\ast}$ then consider the sequence $g=(g_{x})\in\Lambda$ where each $g_{x}$ is either $1$ or $0$, according as $x\in\Su(f)$ or $x\notin\Su(f)$. Then clearly $\Su(g)=\Su(f)\in M$ and so $g\in M^{\ast}$. Moreover $f(1-g)=0$. This yields that $f\in\mathfrak{p}$. Hence, $M^{\ast}$ is a minimal prime of $\Lambda$. The map $\phi$ is continuous, since $\phi^{-1}\big(\Min(\Lambda)\cap D(f)\big)=D\big(\Su(f)\big)$ for all $f\in\Lambda$. Clearly $\mathfrak{m}^{\ast}_{x}=\mathfrak{p}_{x}$ for all $x\in X$, where $\mathfrak{m}_{x}=\mathcal{P}(X\setminus\{x\})$ and for $\mathfrak{p}_{x}$ see the above of Theorem \ref{Theorem STIII}. This shows that $\eta=\phi\circ\eta'$ where $\eta:X\rightarrow\Min(\Lambda)$ and $\eta':X\rightarrow\Spec\mathcal{P}(X)$ are the canonical maps. Therefore, by the universal property of the Stone-\v{C}ech compactification and by taking into account Theorem \ref{Theorem STIII} and Corollary \ref{Corollary STII}, we deduce that $\phi$ is a homeomorphism.
\end{proof}

\begin{theorem}\label{Theorem Ultra 2} The map $\psi:\Spec\mathcal{P}(X)\rightarrow\Max(\Gamma)$ given by $M\rightsquigarrow M^{\flat}$ is a homeomorphism.
\end{theorem}

\begin{proof} By the proof of Theorem \ref{Theorem IUltra}(iii), $M^{\flat}$ is a maximal ideal of $\Gamma$. Hence, the above map is well-defined. It is also continuous, since $\psi^{-1}\big(\Max(\Gamma)\cap D(f)\big)=D\big(\Omega(f)\big)$ for all $f\in\Gamma$. Moreover $\mathfrak{m}^{\flat}_{x}=\mathfrak{M}_{x}$ for all $x\in X$, where $\mathfrak{m}_{x}=\mathcal{P}(X\setminus\{x\})$ and for $\mathfrak{M}_{x}$ see Theorem \ref{Theorem STI}. This shows that
$\eta=\psi\circ\eta'$ where $\eta:X\rightarrow\Max(\Gamma)$ and $\eta':X\rightarrow\Spec\mathcal{P}(X)$ are the canonical maps. Thus, by the universal property of the Stone-\v{C}ech compactification and by taking into account Corollary \ref{Corollary STII} and Theorem \ref{Theorem STI}, we deduce that $\psi$ is a homeomorphism.
\end{proof}

\section{Alexandroff compactification}

Let $\mathcal{R}$ be the set of all subsets of a set $X$ which are either finite or cofinite (i.e. its complement is finite). Then clearly $\mathcal{R}$ is a subring of the power set ring $\mathcal{P}(X)$. Recall that if $x\in X$ then $\mathfrak{m}_{x}=\mathcal{P}(X\setminus\{x\})$ is a maximal ideal of $\mathcal{P}(X)$. In the following result the maximal ideals of $\mathcal{R}$ are characterized.

\begin{theorem}\label{Theorem DIV} Let $X$ be an infinite set. Then the maximal ideals of $\mathcal{R}$ are precisely $\Fin(X)$ or of the form $\mathfrak{m}_{x}\cap\mathcal{R}$ where $x\in X$.
\end{theorem}

\begin{proof} First we have to show that $\Fin(X)$ is a maximal ideal of $\mathcal{R}$. Clearly $\Fin(X)\neq\mathcal{R}$ since $X$ is infinite. If there exists an ideal $I$ of $\mathcal{R}$ strictly containing $\Fin(X)$ then we may choose some $A\in I$ such that $A\notin\Fin(X)$. It follows that $A^{c}\in\Fin(X)$ and so $1=A+A^{c}\in I$. Hence, $\Fin(X)$ is a maximal ideal of $\mathcal{R}$. Conversely, let $M$ be a maximal ideal of $\mathcal{R}$ such that $M\neq\mathfrak{m}_{x}\cap\mathcal{R}$ for all $x\in X$. It follows that $A_{x}:=X\setminus\{x\}\in\mathcal{R}\setminus M$ for all $x\in X$. But $\{x\}.A_{x}=0\in M$. Therefore $\{x\}\in M$ for all $x\in X$. This yields that $\Fin(X)\subseteq M$ and so $\Fin(X)=M$.
\end{proof}

\begin{remark} Let $X$ be an infinite set. Here we give a second proof to show that $\Fin(X)$ is a maximal ideal of $\mathcal{R}$. There exists a maximal ideal $M$ of $\mathcal{P}(X)$ such that $\Fin(X)\subseteq M$ since $\Fin(X)\neq\mathcal{P}(X)$. We have then $\Fin(X)=M\cap\mathcal{R}$. Thus $\Fin(X)$ is a maximal ideal of $\mathcal{R}$.
\end{remark}

\begin{corollary}\label{Corollary IV} If $X$ is an infinite set, then $\Spec(\mathcal{R})$ is the Alexandroff compactification of the discrete space $X$.
\end{corollary}

\begin{proof} The space $\Spec(\mathcal{R})$ is compact. The map $\eta:X\rightarrow\Spec(\mathcal{R})$ given by $x\rightsquigarrow\mathfrak{m}_{x}\cap\mathcal{R}$ is an open embedding. Because by Theorem \ref{Theorem DIV}, $\D(\{x\})=\{\mathfrak{m}_{x}\cap\mathcal{R}\}$ for all $x\in X$.
Now if $A$ is a subset of $X$
then $\eta(A)=\bigcup\limits_{x\in A}\D(\{x\})$. If $U$ is an open neighborhood of $\Fin(X)$ in $\Spec(\mathcal{R})$ then $U^{c}$ is a finite set. Hence, $\eta(X)$ is a dense subspace of $\Spec(\mathcal{R})$.
\end{proof}

Clearly $\Fin(X)$ is the point at infinity of the (infinite) discrete space $X$. \\

Note that, unlike the Stone-\v{C}ech compactification which exists for any topological space, the Alexandroff compactification does not necessarily exist for any space.

\section{Totally disconnected compactifications}

In this section it is shown that every totally disconnected compactification of a discrete space $X$ is precisely of the form $\Spec(\mathcal{R}')$ where the ring $\mathcal{R}'$ satisfies in the extensions of rings $\mathcal{R}\subseteq\mathcal{R}'\subseteq\mathcal{P}(X)$. For $\mathcal{R}$ see \S7.

\begin{lemma}\label{Lemma III} Let $f:X\rightarrow Y$ be a continuous map of topological spaces such that $f(X)$ is a dense subspace of $Y$. Then the induced map $\Clop(f):\Clop(Y)\rightarrow\Clop(X)$ is injective.
\end{lemma}

\begin{proof} Let $A$ be a clopen of $Y$ such that $f^{-1}(A)=\emptyset$. If $A$ is non-empty then $A\cap f(X)$ is non-empty. But this is a contradiction. \\
As a second proof, let $D_{1}$ and $D_{2}$ be two clopens of $Y$ such that $f^{-1}(D_{1})=f^{-1}(D_{2})$. Suppose there exists some $y\in D_{1}$ such that $y\notin D_{2}$. It follows that $(D_{1}\cap D_{2}^{c})\cap f(X)\neq\emptyset$. Hence there exists some $x\in X$ such that $f(x)\in D_{1}\cap D_{2}^{c}$. But this is a contradiction. Therefore $D_{1}=D_{2}$.
\end{proof}

\begin{theorem}\label{Theorem IV} Every totally disconnected  compactification of a discrete space $X$ is precisely of the form $\Spec(\mathcal{R}')$ where the ring $\mathcal{R}'$ satisfies in the extensions of rings $\mathcal{R}\subseteq\mathcal{R}'\subseteq\mathcal{P}(X)$.
\end{theorem}

\begin{proof} It is easy to see that for any such ring $\mathcal{R}'$ then $\Spec(\mathcal{R}')$ together with the canonical open embedding  $\eta:X\rightarrow\Spec(\mathcal{R}')$ which sends each point $x\in X$ into $\mathfrak{m}_{x}\cap\mathcal{R}'$ is a totally disconnected compactification of the discrete space $X$. Conversely, let $(\widetilde{X},\eta)$ be a totally disconnected compactification of a discrete space $X$. By \cite[Corollary 5.4]{Tarizadeh}, the space $\widetilde{X}$ is homeomorphic to $\Spec(R)$ where $R=\Clop(\widetilde{X})$. By Lemma \ref{Lemma III}, the induced map $\Clop(\eta):R
\rightarrow\Clop(X)=\mathcal{P}(X)$ is an injective ring map. So the ring $R$ is isomorphic to $\mathcal{R}'$, the image of $\Clop(\eta)$. It remains to show that $\mathcal{R}\subseteq\mathcal{R}'$. Take $A\in\mathcal{R}$. If $A$ is finite then $D:=\eta(A)=\bigcup\limits_{x\in A}\{\eta(x)\}$ is a closed subset of $\widetilde{X}$ and so $D\in\Clop(\widetilde{X})$. Therefore $A=\eta^{-1}(D)\in\mathcal{R}'$. But if $A$ is cofinite then the above argument shows that $A^{c}\in\mathcal{R}'$, and so $A=1-A^{c}\in\mathcal{R}'$.
\end{proof}

\begin{remark} If $(\widetilde{X},\eta)$ is an arbitrary compactification of a discrete space $X$ then  by \cite[Theorem 5.2]{Tarizadeh}, the space of connected components $\pi_{0}(\widetilde{X})$ is homeomorphic to $\Spec(\mathcal{R}')$ where $\mathcal{R}'=\Clop(\widetilde{X})$. Also $\mathcal{R}'$, via the ring map $\Clop(\eta)$, can be viewed as a subring of $\mathcal{P}(X)$ and containing $\mathcal{R}$.
Note that there are compactifications of a discrete space which are not totally disconnected.
\end{remark}

\section{Semigroup structure on $\beta X$}

In this section, $\beta X=\Spec\mathcal{P}(X)$ together with the canonical map $\eta:X\rightarrow\beta X$ denotes the Stone-\v{C}ech compactification of the discrete space $X$. If $f:X\rightarrow Y$ is a function then by Corollary \ref{Corollary STII}, there exists a unique continuous function $\beta f:\beta X\rightarrow\beta Y$ such that $(\beta f)(\mathfrak{m}_{x})=\mathfrak{m}_{f(x)}$ for all $x\in X$. This yields that $\beta f=\mathcal{P}(f)^{\ast}$. In particular, if $f:X\rightarrow Y$ is injective then $\beta f:\beta X\rightarrow\beta Y$ is as well. \\

Let $(S,\ast)$ be a semigroup such that $S$ is a topological space. Then equip $S\times S$ with the product topology. If the operation $\ast:S\times S\rightarrow S$ is continuous, then $(S,\ast)$ is called a topological semigroup. If the operation $\ast$ is not continuous, then this leads us to a weaker notion. Indeed, the pair $(S,\ast)$ is called a \emph{left topological semigroup} if the operation $\ast$ is left semi-continuous. That is, for each $p\in S$ then the map $\ell_{p}:S\rightarrow S$ given by $x\rightsquigarrow p\ast x$ is continuous. The right topological semigroup is defined dually. Obviously every topological semigroup is both right topological and left topological semigroup. We have then the following interesting result.

\begin{theorem}\label{Theorem DV} The operation of every commutative semigroup $(X,.)$ can be extended uniquely to an operation $\ast$ on $\beta X$ such that: $(\beta X,\ast)$ is a left topological semigroup, the canonical map $\eta:X\rightarrow\beta X$ is a morphism of semigroups and $\mathfrak{m}_{x}\ast M=M\ast\mathfrak{m}_{x}$ for all $M\in\beta X$ and $x\in X$. If moreover $e$ is the identity of $X$, then $\mathfrak{m}_{e}$ is the identity of $\beta X$.
\end{theorem}

\begin{proof} If $x\in X$ then by Theorem \ref{Theorem DI}, there exists a unique continuous function $\phi_{x}:\beta X\rightarrow\beta X$ such that $\phi_{x}(\mathfrak{m}_{y})=\mathfrak{m}_{x.y}$ for all $y\in X$. For a fixed $M\in\beta X$, again by Theorem \ref{Theorem DI}, there exists a unique continuous map $\theta_{M}:\beta X\rightarrow\beta X$ such that $\theta_{M}(\mathfrak{m}_{x})=\phi_{x}(M)$ for all $x\in X$. Now we define the operation $\ast$ on $\beta X$ as $M\ast N=\theta_{M}(N)$.
Then we show this operation is associative. To prove this it suffices to show that $\theta_{M}\circ\theta_{N}=\theta_{L}$ for every $M,N\in\beta X$ with $L=\theta_{M}(N)$. To see this it will be enough to show that $\theta_{M}\circ\phi_{x}=\phi_{x}\circ\theta_{M}$ for all $M\in\beta X$ and $x\in X$. But to see the latter it suffices to show that $\theta_{M}\circ\phi_{x}$ and $\phi_{x}\circ\theta_{M}$ agree on $\eta(X)$, (recall that if two continuous maps into a Hausdorff space agree on a dense subspace of the domain, they are equal). This reduces to show that $\phi_{x}\circ\phi_{y}=\phi_{x.y}$ for all $x,y\in X$. Finally, to see this it suffices to show that $(\phi_{x}\circ\phi_{y})(\mathfrak{m}_{z})=\phi_{xy}(\mathfrak{m}_{z})$ for all $z\in X$. But the latter obviously holds since the operation of $X$ is associative. Clearly $\ell_{M}=\theta_{M}$ for all $M\in\beta X$. Hence, $(\beta X,\ast)$ is a left topological semigroup. The map $\eta$ is a morphism of semigroups since $\phi_{x}=\theta_{\mathfrak{m}_{x}}$ for all $x\in X$. This also yields that $\mathfrak{m}_{x}\ast M=M\ast\mathfrak{m}_{x}$ for all $M\in\beta X$ and $x\in X$. To see the uniqueness of $\ast$, suppose there is another operation $\ast'$ on $\beta X$ such that $(\beta X,\ast')$ is a left topological semigroup, the canonical map $\eta:X\rightarrow(\beta X,\ast')$ is a morphism of semigroups and $\mathfrak{m}_{x}\ast'M=M\ast'\mathfrak{m}_{x}$ for all $M\in\beta X$ and $x\in X$. Then clearly for each $x\in X$, the maps $\ell_{\mathfrak{m}_{x}}$ and $\ell'_{\mathfrak{m}_{x}}$ agree on $\eta(X)$, hence they are equal. It follows that for each $M\in\beta X$, then $\ell_{M}$ and $\ell'_{M}$ agree on $\eta(X)$, hence they are equal. The latter implies that $\ast=\ast'$. Finally, if $e$ is the identity element of $X$ then $\phi_{e}$ is the identity map. It follows that $\mathfrak{m}_{e}$ is the identity element of $\beta X$.
\end{proof}

Note that the operation $\ast$ of Theorem \ref{Theorem DV} is not necessarily commutative. Hence, we may define a new operation on $\beta X$ as $M\ast'N:=\theta_{N}(M)=N\ast M$. Then it is easy to see that $(\beta X,\ast')$ is a right topological semigroup. Therefore we may consider $\beta X$ as left topological or right topological semigroup, depending on the preferred construction, but never both (specially when $X$ is an infinite set). \\

In the proof of Theorem \ref{Theorem DV}, we have $\phi_{x}=\mathcal{P}(f_{x})^{\ast}$ for all $x\in X$ where the function  $f_{x}:X\rightarrow X$ is defined by $f_{x}(y)=x.y$. By \cite[Theorem 5.6 ]{Tarizadeh}, there exists a (unique) morphism of rings $h_{M}:\mathcal{P}(X)\rightarrow\mathcal{P}(X)$ such that $\theta_{M}=\Spec(h_{M})$. In the following result, the rule of this morphism is determined explicitly.

\begin{theorem}\label{Theorem the4ta} Let $(X,.)$ be a commutative semigroup and $M$ a maximal ideal of $\mathcal{P}(X)$. Then the map $\zeta_{M}:\mathcal{P}(X)\rightarrow\mathcal{P}(X)$ given by
$A\rightsquigarrow\{x\in X: f^{-1}_{x}(A)\notin M\}$ is a morphism of rings and $\theta_{M}=\Spec(\zeta_{M})$.
\end{theorem}

\begin{proof} It is not hard to see that the map $\zeta_{M}$ is actually a morphism of rings. To see $\theta_{M}=\Spec(\zeta_{M})$ it suffices to show that $\theta_{M}(\mathfrak{m}_{x})=\zeta^{-1}_{M}(\mathfrak{m}_{x})$ for all $x\in X$.
\end{proof}

By the category of left topological monoids we mean a category whose objects are the left topological monoids and whose morphisms are the continuous morphisms of monoids.

\begin{theorem}\label{Proposition Elman-Asena} The assignments $X\rightsquigarrow\beta X$ and $h\rightsquigarrow\beta f$ form a faithful covariant functor from the category of commutative monoids to the category of left topological monoids.
\end{theorem}

\begin{proof} By the universal property of the Stone-\v{C}ech compactification, it is a functor provided that we could prove that if $f:X\rightarrow Y$ is a morphism of commutative monoids then $\beta f:\beta X\rightarrow\beta Y$ is a morphism of monoids. Clearly $\beta f$ preserves the identities. It remains to show that $(\beta f)\circ\theta_{M}=\theta_{M'}\circ(\beta f)$ for all $M\in\Spec\mathcal{P}(X)$ with $M'=(\beta f)(M)$, for the notations see the proof of Theorem \ref{Theorem DV}. To see this it suffices to show that these functions agree on $\eta(X)$. To see the latter it will be enough to show that  $(\beta f)\circ\phi_{x}=\phi_{f(x)}\circ(\beta f)$ for all $x\in X$. Clearly these maps agree on $\eta(X)$, hence they are equal.
\end{proof}

\section{Absolutely flatness of the total ring of fractions}

Theorems \ref{Theorem Asena-Elman} and \ref{Proposition Quentel} provide new and simple proofs to the main results of  \cite[Theorem 2.9]{Glaz 2}, \cite[Chap I, Theorem 4.5]{Huckaba}, \cite[Proposition 1.4]{Matlis} and \cite[Proposition 9]{Quentel}. In the following results, $T(R)$ denotes the total ring of fractions of a ring $R$.

\begin{theorem}\label{Theorem Asena-Elman} Let $R$ be a ring. Then $T(R)$ is absolutely flat if and only if the following two conditions hold. \\
$\mathbf{(i)}$ $R$ is reduced and $\Min(R)$ is Zariski compact. \\
$\mathbf{(ii)}$ Every finitely generated and faithful ideal of $R$ contains a non zero-divisor of $R$.
\end{theorem}

\begin{proof} Assume $T(R)$ is absolutely flat. Then $R$ is reduced, since every absolutely flat ring and so each subring are reduced. By \cite[Lemma 3.4]{A-Tarizadeh-Racsam}, $\Min(R)$ is Zariski compact.
If $I=(f_{1},...,f_{n})$ is a finitely generated and faithful ideal of $R$ then for each $i$, there exists a non zero-divisor $g_{i}$ of $R$ such that $f_{i}(g_{i}-f_{i}h_{i})=0$ for some $h_{i}\in R$. It follows that $(g_{1}-f_{1}h_{1})...(g_{n}-f_{n}h_{n})=0$ and so $g_{1}...g_{n}\in I$. Conversely, if $f\in R$ then it will be enough to find a non zero-divisor $g$ of $R$ such that $fg=f^{2}h$ for some $h\in R$. Setting $X=\{\mathfrak{p}\in\Min(R): f\in\mathfrak{p}\}$. If $\mathfrak{p}\in X$ then there exists some $x_{\mathfrak{p}}\in R\setminus\mathfrak{p}$ such that $fx_{\mathfrak{p}}=0$. It follows that $\Min(R)\subseteq D(f)\cup\big(\bigcup\limits_{\mathfrak{p}\in X}D(x_{\mathfrak{p}})\big)$. Using the quasi-compactness of $\Min(R)$, then we may write $\Min(R)\subseteq D(f)\cup\big(\bigcup\limits_{i=1}^{n}D(x_{i})\big)$ and that $fx_{i}=0$ for all $i$. Therefore $I=(f,x_{1},...,x_{n})$ is a faithful ideal of $R$, because suppose $rI=0$, if $\mathfrak{p}\in\Min(R)$ then $r\in\mathfrak{p}$ and so $r\in\bigcap\limits_{\mathfrak{p}\in\Min(R)}\mathfrak{p}=
\sqrt{0}=0$. Hence, $I$ contains a non zero-divisor $g$ of $R$. Thus we may write $g=fh+\sum\limits_{i=1}^{n}r_{i}x_{i}$ where $h,r_{1},...,r_{n}\in R$. This yields that $fg=f^{2}h$.
\end{proof}

\begin{remark} Let $R$ be a ring. It is easy to see that if at least one of the coefficients of a polynomial $f\in R[x]$ is a non zero-divisor of $R$, then $f$ is a non zero-divisor of $R[x]$. But the converse does not hold. As an example, take $R=\mathbb{Z}/6\mathbb{Z}$ then $f=2+3x$ is a non zero-divisor of $R[x]$, but all of its coefficients are zero-divisors of $R$. This observation shows that if $T(R[x])$ is zero dimensional (or, an absolutely flat ring) then the same assertion
does not necessarily hold for $T(R)$.
\end{remark}

\begin{corollary} Let $R$ be a reduced ring such that $\Min(R)$ is a finite set. Then $T(R)$ is absolutely flat.
\end{corollary}

\begin{proof} Let $I=(f_{1},...,f_{n})$ be a faithful ideal of $R$ and setting $S:=R\setminus Z(R)$. If $I\cap S=\emptyset$ then there exists a prime ideal $\mathfrak{p}$ of $R$ such that $I\subseteq\mathfrak{p}$ and $\mathfrak{p}\cap S=\emptyset$. It follows that $\mathfrak{p}\subseteq Z(R)=\bigcup\limits_{\mathfrak{q}\in\Min(R)}\mathfrak{q}$. Thus by the Prime Avoidance Lemma (cf. \cite[Theorem 2.2]{Tarizadeh-Chen}), $\mathfrak{p}\in\Min(R)$. So for each $i$, there exists some $g_{i}\in R\setminus\mathfrak{p}$ such that $f_{i}g_{i}=0$. Therefore $gI=0$ where $g=g_{1}\ldots g_{n}$. But this is a contradiction. Hence, $I$ admits a non zero-divisor of $R$. Therefore by Theorem \ref{Theorem Asena-Elman}, $T(R)$ is absolutely flat.
\end{proof}

The following two results are easily deduced from the above corollary.

\begin{corollary} Let $R$ be a reduced ring such that $\Spec(R)$ is a noetherian space with respect to the Zariski topology. Then $T(R)$ is absolutely flat.
\end{corollary}

\begin{corollary} Let $R$ be a reduced and noetherian ring. Then $T(R)$ is absolutely flat.
\end{corollary}

\begin{corollary} Let $R$ be a ring. Then $T(R[x])$ is absolutely flat if and only if $R$ is reduced and $\Min(R)$ is Zariski compact.
\end{corollary}

\begin{proof} It is interesting to notice that for any ring $R$, then every finitely generated and faithful ideal of $R[x]$ contains a non zero-divisor of $R[x]$. Then apply Theorem \ref{Theorem Asena-Elman}.
\end{proof}

\begin{corollary} Let $R$ be a ring. If $T(R)$ is absolutely flat, then $T(R[x])$ is as well.
\end{corollary}

\begin{remark} Here we prove the reverse implication of Theorem \ref{Theorem Asena-Elman} by an alternative approach. Clearly a ring $R$ is reduced if and only if $T(R)$ is reduced. To prove that $T(R)$ is zero dimensional it will be enough to show that every prime ideal $\mathfrak{q}$ of $R$ which does not meet $R\setminus Z(R)$, then it is a minimal prime. To see this it suffices to show that $\mathfrak{q}\subseteq\mathfrak{p}$ for some $\mathfrak{p}\in\Min(R)$. Suppose  $\mathfrak{q}\nsubseteq\mathfrak{p}$ for all $\mathfrak{p}\in\Min(R)$, then there exists some $x_{\mathfrak{p}}\in\mathfrak{q}$ such that $x_{\mathfrak{p}}\notin\mathfrak{p}$. So $\Min(R)\subseteq
\bigcup\limits_{\mathfrak{p}\in\Min(R)}D(x_{\mathfrak{p}})$. By the quasi-compactness of the minimal spectrum, we may write $\Min(R)\subseteq
\bigcup\limits_{i=1}^{n}D(x_{i})$ where $x_{i}:=x_{\mathfrak{p}_{i}}$ for all $i$. Then $I=(x_{1},\ldots,x_{n})\subseteq\mathfrak{q}\subseteq Z(R)$. So $I$ is not a faithful ideal, i.e., $\Ann(I)\neq0$. Thus we may choose some nonzero $a$ in $\Ann(I)$. Now if $\mathfrak{p}\in\Min(R)$ then $x_{i}\notin\mathfrak{p}$ for some $i$. But $ax_{i}=0$. Thus $a\in\mathfrak{p}$. So $a\in\bigcap\limits_{\mathfrak{p}\in\Min(R)}\mathfrak{p}=0$ which is a contradiction.
\end{remark}

\begin{theorem}\label{Proposition Quentel} For a reduced ring $R$ the following assertions are equivalent. \\
$\mathbf{(i)}$ $T(R)$ is absolutely flat. \\
$\mathbf{(ii)}$ If $f\in R$ there exists some $g\in R$ such that $fg=\Ann(Rf+Rg)=0$. \\
$\mathbf{(iii)}$ If an ideal $I$ of $R$ is contained in $Z(R)$, then $I\subseteq\mathfrak{p}$ for some $\mathfrak{p}\in\Min(R)$.
\end{theorem}

\begin{proof} $\mathbf{(i)}\Rightarrow\mathbf{(ii)}:$ There exists a non zero-divisor $s\in R$ such that $fs=f^{2}h$ for some $h\in R$. Then $g:=fh-s$ is the desired element. \\ $\mathbf{(ii)}\Rightarrow\mathbf{(i)}:$ It suffices to show that $h:=f-g$ is a non zero-divisor of $R$. Suppose $rh=0$ and $\mathfrak{p}$ is a minimal prime ideal of $R$ such that $r\notin\mathfrak{p}$. Then $f,g\in\mathfrak{p}$ and so there exist $f',g'\in R\setminus\mathfrak{p}$ such that $ff'=gg'=0$. This yields that $f'g'\in\Ann(Rf+Rg)=0$ which is a contradiction. Hence, $r\in\bigcap\limits_{\mathfrak{p}\in\Min(R)}\mathfrak{p}
=\sqrt{0}=0$. \\
$\mathbf{(i)}\Rightarrow\mathbf{(iii)}:$ Suppose $I\nsubseteq\mathfrak{p}$ for all $\mathfrak{p}\in\Min(R)$, then we may choose some $x_{\mathfrak{p}}\in I\setminus\mathfrak{p}$. Using the quasi-compactness of $\Min(R)$, then we may write $\Min(R)\subseteq\bigcup\limits_{i=1}^{n}D(x_{i})$ where $x_{i}\in I$ for all $i$. It follows that $J=(x_{1},...,x_{n})$ is a faithful ideal of $R$. Thus by Theorem \ref{Theorem Asena-Elman}, $J$ admits a non zero-divisor which is a contradiction. \\
$\mathbf{(iii)}\Rightarrow\mathbf{(i)}:$ Suppose $\Min(R)\subseteq\bigcup\limits_{i\in S}D(f_{i})$ where $f_{i}\in R$ for all $i$. Then by the hypothesis, the ideal $(f_{i}: i\in S)$ admits a non zero-divisor $g$ of $R$. So there exists a finite subset $S'$ of $S$ such that  $g=\sum\limits_{i\in S'}r_{i}f_{i}$ where $r_{i}\in R$ for all $i\in S'$. This yields that $\Min(R)\subseteq\bigcup\limits_{i\in S'}D(f_{i})$, since otherwise we may find some $\mathfrak{p}\in\Min(R)$ such that $g\in \mathfrak{p}$, but this is impossible since $\mathfrak{p}\subseteq Z(R)$. Hence, $\Min(R)$ is quasi-compact. Now let $I$ be a finitely generated and  faithful ideal of $R$. If $I\subseteq Z(R)$ then $I$ is contained in a minimal prime ideal $\mathfrak{p}$ of $R$. Thus we may find some $s\in R\setminus\mathfrak{p}$ such that $sI=0$, which is a contradiction. So $I$ admits a non zero-divisor of $R$. Therefore by Theorem \ref{Theorem Asena-Elman}, $T(R)$ is absolutely flat.
\end{proof}

\begin{proposition} Let $R$ be a ring. If $T(R[x])$ is a zero dimensional ring, then $\Min(R)$ is Zariski compact.
\end{proposition}

\begin{proof} For any ring $R$, the minimal prime ideals of $R[x]$ are precisely of the form $\mathfrak{p}[x]$ where $\mathfrak{p}$ is a minimal prime ideal of $R$. Hence, the map $\phi:\Min(R)\rightarrow\Min(R[x])$ given by $\mathfrak{p}\rightsquigarrow\mathfrak{p}[x]$ is bijective.
This map is continuous, because if $f=\sum\limits_{i=0}^{n}f_{i}x^{i}\in R[x]$ with the $f_{i}\in R$, then
$\phi^{-1}(U)=\bigcup\limits_{i=0}^{n}U_{i}$ where $U=\Min(R[x])\cap D(f)$ and $U_{i}=\Min(R)\cap D(f_{i})$. The converse of $\phi$ is also continuous, because it is induced by the ring extension $R\subseteq R[x]$. Therefore $\phi$ is a homeomorphism.
By \cite[Lemma 3.4]{A-Tarizadeh-Racsam}, $\Min(R[x])$ is Zariski compact. Hence, $\Min(R)$ is Zariski compact.
\end{proof}

\end{document}